\documentclass{elsarticle}

\usepackage[final,                     
color,notref,notcite]{showkeys}         
\usepackage{amsmath}
\usepackage{amsfonts} 
\usepackage{bm}
\usepackage{tikz}
\usepackage{pgfplots}
\usepgflibrary{arrows}
\usepackage{rotating} 
\graphicspath{{figs/}}

\newtheorem{thm}{Theorem}
\newtheorem{lem}[thm]{Lemma}
\newproof{pf}{Proof}

\newcommand{\grad}{\mathop{\rm grad}\nolimits}
\renewcommand{\div}{\mathop{\rm div}\nolimits}
\newcommand{\const}{\mathop{\rm const}\nolimits}

\journal{arXiv} 

\begin{document}

\begin{frontmatter}

\title{Splitting Schemes for Some Second-Order Evolution Equations}

\author[nsi,uni]{Petr N. Vabishchevich\corref{cor}}
\ead{vabishchevich@gmail.com}

\address[nsi]{Nuclear Safety Institute, Russian Academy of Sciences, Moscow, Russia}
\address[uni]{North-Eastern Federal University, Yakutsk, Russia}

\cortext[cor]{Corresponding author}

\begin{abstract}

We consider the Cauchy problem for a second-order evolution equation,
in which the problem operator is the sum of two self-adjoint operators.
The main feature of the problem is that one of the operators is represented in the form
of the product of operator $A$ by its conjugate $A^*$.
Time approximations are carried out so that
the transition to a new level in time was associated with a separate solution of problems for
operators $A$ and $A^*$, not their products.
The construction of unconditionally stable schemes is based on general results of
the theory of stability (correctness) of operator-difference schemes in Hilbert spaces
and is associated with the multiplicative perturbation of the problem operators, which lead to stable implicit schemes.
As an example, the problem of the dynamics of a thin plate on an elastic foundation is considered.

\end{abstract}


\end{frontmatter}

\section{Introduction}

Many applied problems lead to the need for an approximate solution of the Cauchy problem
for second-order evolution equations.
As a typical example, we note the dynamic
problems of solid mechanics \cite{fung2017classical}.
A class of problems can be distinguished, a characteristic feature of which is that the main part of the problem operator
is the product of two operators. For example, when considering models of thin plates
we have a biharmonic operator, the product of two Laplace operators.

Unconditionally stable schemes for these problems are built based on implicit approximations in time \cite{SamarskiiTheory,LeVeque2007}.
In the theory of stability (correctness) of operator-difference schemes \cite{SamarskiiGulin1973,SamarskiiMatusVabischevich2002}
the most complete results were obtained on the stability of two-level and three-level schemes in Hilbert spaces.
The computational complexity of solving the Cauchy problem on a new level in time using implicit schemes may be unacceptable.
Therefore, various approaches are being developed to obtain computationally simpler problems when solving non-stationary problems.

Simplification of the problem on a new level is often implemented for evolutionary problems when the problem operator is represented in the form
the sums are more simple. For such problems, additive operator-difference schemes are constructed,
which are related to one or another inhomogeneous approximation in time for individual operator terms.
The traditional approach is based on explicit-implicit approximations (IMEX methods) \cite{ascher1995implicit,HundsdorferVerwer2003}
when one part of the problem operator is taken from the lower level in time (explicit approximation), and the other --- from the upper one (implicit approximation).
This idea of time approximation is implemented most consistently when constructing splitting schemes \cite{Marchuk1990,VabishchevichAdditive}.
In this case, the transition to a new level in time is carried out by solving evolutionary problems for individual operator terms.

One more class of evolutionary problems can also be noted, in which the problem operator is represented as a product of
two or more operators. An example is nonstationary problems with a variable weighting factor, the study of which is
held in \cite{SamarskiiMatusVabischevich2002,SamVabGul}.
Special time approximations are constructed to simplify the problem on a new time level.
For example, paper \cite{VabFact} builds schemes that are based
on the solution of a discrete problem on a new time level with one operator factor.

In this paper, we consider the Cauchy problem for a second-order evolution equation in which the problem operator includes
the product of operator $A$ by its conjugate $A^*$.
Unconditionally stable schemes are constructed based on a perturbation of both the $A$ operator and the $A^*$ operator.
In this case, the computational implementation is associated with the separate solution of problems for
operators $A$ and $A^*$, not their products.

The article is organized as follows.
Statement of the Cauchy problem for a second-order evolution equation, which includes
the product of the operator's $A$ and $A^*$ is given in Section 2.
Section 3 describes a general approach to constructing unconditionally stable schemes for second-order evolution equations
based on multiplicative perturbation of the operator(s) of the problem.
Splitting schemes for the evolutionary problem, when the problem operator includes $A^*A$, are constructed in Section 4.
In Section 5, the results obtained are applied to the model problem of the dynamics of a thin plate on an elastic foundation.
The results of our work are summarized in Section 6.

\section{Problem statement}

The Cauchy problem for a second-order evolution equation is considered in a finite-dimensional Hilbert space
$H$. In order not to clutter up the presentation with technical details, we restrict ourselves to a homogeneous equation when
\begin{equation}\label{2.1}
 \frac{d^2 w}{d t^2} + A^* A w + B w = 0,
 \quad 0 < t \leq T, 
\end{equation} 
\begin{equation}\label{2.2}
 w(0)= w^0 ,
 \quad \frac{d w}{d t} (0) = \widetilde{w}^0 .
\end{equation}
We will assume that the operator's $A$ and $B$ in (\ref{2.1}) are constant (do not depend on $t$),
and operator $B$ is self-adjoint and non-negative:
\begin{equation}\label{2.3}
 B = B^* \geq  0 . 
\end{equation} 

We arrive at the problem (\ref{2.1})--(\ref{2.3}), for example, after discretization by
spatial variables in the numerical solution of initial boundary value problems
for hyperbolic equations. The key feature of the problem under consideration is associated with  operator $A$,
so that it enters the equation (\ref{2.1}) as the product $A^* A$.
An example of such a construction is the biharmonic operator ($A = A^*$).

The scalar product for $u, v \in H$ is $(u, v)$, and the norm is $\| u \| = (u, u)^{1/2}$.
Let us define a Hilbert space $H_S$ with scalar product and norm $(u, v)_S = (S u, v), \ \|u\|_S = (u, v)_S^{1/2}$,
which is generated by the self-adjoint and positive operator $S$.

The subject of our consideration is time approximation for equation (\ref{2.1}).
We focus on unconditionally stable schemes for an approximate solution to the problem (\ref{2.1})--(\ref{2.3}),
which are convenient for computational implementation. When obtaining the corresponding stability estimates
we compare them with a priori estimates for the differential problem.

We multiply the equation (\ref{2.1}) scalarly in $H$ by $d w/dt$ and obtain
\[
 \frac{d}{d t} \left ( \left \| \frac{d w}{d t} \right \|^2 + \|A w\|^2 + \|w\|_B^2 \right ) = 0 . 
\] 
This equality implies the estimate
\begin{equation}\label{2.4}
 \left \| \frac{d w}{d t}(t) \right \|^2 + \|A w(t)\|^2 + \|w(t)\|_B^2 = \|\widetilde{w}^0\|^2 + \|A w^0\|^2 + \|w^0\|_B^2 ,
\end{equation} 
which ensures stability with respect to the initial data of the solution to the problem (\ref{2.1})--(\ref{2.3}).

We will use a uniform, for simplicity, grid in time with
step $\tau$ and notation $u^n = u(t^n), \ t^n = n \tau$,
$n = 0,\ldots,N, \ N\tau = T$.
As a basic scheme for the numerical solution of the problem (\ref{2.1})--(\ref{2.3})
we will use a three-level scheme with weights ($\sigma = \mathrm{const}$):
\begin{equation}\label{2.5}
\begin{split}
 \frac{u^{n+1} - 2 u^{n} + u^{n-1}}{\tau^2 } + & (A^* A + B) (\sigma u^{n+1} + (1-2\sigma) u^{n} + \sigma u^{n-1}) = 0, \\
 & \quad n = 1,\ldots,N-1 , 
\end{split}
\end{equation} 
when setting the initial conditions
\begin{equation}\label{2.6}
 u^0 = w^0 , 
 \quad u^1 = \overline{w}^1 . 
\end{equation} 
For the second initial condition on the solutions of the problem (\ref{2.1}), (\ref{2.2}) we put
\[
 \left (I + \frac{\tau^2}{2}  (A^* A + B) \right )  \overline{w}^1 = w^0 + \tau \widetilde{w}^0. 
\] 
Difference scheme (\ref{2.5}), (\ref{2.6}) approximates (\ref{2.1}), (\ref{2.2})
with second-order in $\tau$.

Our consideration is based on the use of general results in
the theory of stability (correctness) of operator-difference schemes
in Hilbert spaces \cite{SamarskiiGulin1973,SamarskiiMatusVabischevich2002}.
The main statement on the stability of three-level schemes for the problems under consideration is formulated as follows. 

\begin{lem}\label{l-1}
Let in a three-level scheme
\begin{equation}\label{2.7}
 C \frac{u^{n+1} - 2 u^n + u^{n-1}}{\tau^2} + D u^{n} = 0, 
 \quad n = 1,\ldots, N-1 , 
\end{equation} 
when specifying (\ref{2.6}), the operators
\begin{equation}\label{2.8}
 C = C^* > 0,
 \quad D = D^* > 0 . 
\end{equation} 
Then at
\begin{equation}\label{2.9}
 G = C - \frac{\tau^2} {4} D \geq 0
\end{equation} 
the scheme (\ref{2.6})--(\ref{2.8}) is stable and the solution has an a priori equality
\begin{equation}\label{2.10}
 \left \| \frac{u^{n+1}-u^{n}}{\tau } \right \|^2_{G} +
 \left \| \frac{u^{n+1}+u^{n}}{2 } \right \|^2_D 
 = \left \| \frac{\overline{w}^1-w^0}{\tau } \right \|^2_{G} +
 \left \| \frac{\overline{w}^1+w^0}{2 } \right \|^2_D ,
\end{equation} 
for all $n = 1,\ldots, N-1$.
\end{lem}

\begin{pf}
Taking into account
\[
 u^{n} = \frac{u^{n+1} + 2 u^n + u^{n-1}}{4} 
 - \frac{\tau^2}{4} \frac{u^{n+1} - 2 u^n + u^{n-1}}{\tau^2} 
\] 
rewrite (\ref{2.7}) in the form
\begin{equation}\label{2.11}
 G\frac{u^{n+1} - 2 u^n + u^{n-1}}{\tau^2} + D \frac{u^{n+1} + 2 u^n + u^{n-1}}{4} = 0 . 
\end{equation} 
Let's introduce new variables
\[
 s^n = \frac{u^{n} + u^{n-1}}{2} ,
 \quad  r^n = \frac{u^{n} - u^{n-1}}{\tau } ,
\] 
and from (\ref{2.11}) we arrive at the equation
\[
  G\frac{r^{n+1} - r^n}{\tau} + D \frac{s^{n+1} + s^n}{2} = 0 . 
\] 
Let's multiply it by
\[
 2 (s^{n+1} - s^n) = \tau (r^{n+1} + r^n) ,
\] 
what gives
\[
 (G r^{n+1}, r^{n+1}) +  (D s^{n+1}, s^{n+1}) = (G r^{n}, r^{n}) +  (D s^{n}, s^{n}) .
\] 
Returning to the original variables, we arrive at the equality being proved (\ref{2.10}). 
\end{pf}

Application of this lemma to a weighted scheme (\ref{2.3}), (\ref{2.5}), (\ref{2.6}) brings us
to the next statement.

\begin{thm}\label{t-2}
Three-level scheme (\ref{2.3}), (\ref{2.5}), (\ref{2.6}) is unconditionally stable at
$\sigma \geq 1/4$. Under these constraints, for an approximate solution of the problem, the a priori equality (\ref{2.10}),
wherein
\[
 G = I + \left (\sigma - \frac{1}{4} \right) \tau^2 D,
 \quad D = A^* A + B ,   
\] 
and $I$ is the identity operator.
\end{thm}

\begin{pf}
We write (\ref{2.5}) in the form (\ref{2.7}) for
\[
 C = I + \sigma \tau^2 D .
\] 
Conditions (\ref{2.8}) for $\sigma \geq 0$ and (\ref{2.3}) are satisfied, and the inequality (\ref{2.9})
results in $\sigma \geq 1/4$ constraints. Thus, all conditions of the lemma~\ref{l-1} are executed.
\end{pf}

When using the scheme (\ref{2.5}), (\ref{2.6}) on a new $n+1$ level, the problem is solved
\[
 (I + \sigma \tau^2 (A^* A + B)) u^{n+1} = \varphi^n
\]
with known right side $\varphi^n$. The computational complexity of this problem may be unacceptable
and therefore it is necessary to simplify the problem on a new level in time by using
special time approximations.
In our case, we want to ensure the transition to a new level in time by solving individual problems for
operators $A$ and $A^*$, avoiding solving a more complex problem with the product of these operators.

\section{Unconditionally stable schemes with multiplicative regularization} 

The principle of regularization of difference schemes provides great opportunities for
constructing difference schemes of a given quality \cite{VabishchevichAdditive,SamarskiiReg}.
Results of the theory of regularization of difference schemes
are used to improve the quality of the difference scheme due to
introducing regularizers into the operators of the original difference scheme.
The regularization principle for constructing unconditionally stable difference schemes
implemented as follows:

\begin{enumerate}
 \item for the problem under consideration, the simplest
   difference scheme (generating difference scheme),
   not possessing the necessary properties,
   that is, the scheme is conditionally stable or even absolutely
   unstable;
 \item the difference scheme is written in a unified (canonical) form, for
which stability conditions are known;
 \item the quality of the difference scheme (its stability)
   improves due to the perturbation of the difference scheme operators.
\end{enumerate} 
Thus, the principle of regularization of difference schemes is based on the use of
already known general stability conditions, which are given by the theory of stability (correctness) of operator-difference schemes.

Consider the model Cauchy problem for the equation
\begin{equation}\label{3.1}
  \frac {d^2 w}{dt^2} + Q w = 0,
  \quad 0 < t \leq T,
\end{equation}
with a constant, self-adjoint, and positive in $H$ linear operator
$Q$. Following the regularization principle, we first choose some
difference scheme for the problem (\ref{2.2}), (\ref{3.1}), from which we will start. As such
a generating scheme, it is natural to consider the simplest explicit scheme
\begin{equation}\label{3.2}
  \frac{u^{n+1} - 2 u^n + u^{n-1}} {\tau^2} + Q u^n = 0,
  \quad n = 1,\ldots,N-1 ,
\end{equation}
with initial conditions (\ref {2.6}).

To use lemma \ref{l-1}, we write the difference scheme (\ref{3.2}) in the form (\ref{2.7}) with the operators
$C = I, \ D = Q$. Taking into account that $Q \leq \|Q\| I$, from
(\ref{2.9}) we get a time step constraint
\[
 \tau \leq \tau_0 = \frac{2}{\|Q\|^{1/2}} .
\] 
for scheme stability (\ref{2.2}), (\ref{3.1}).
 
By (\ref{2.9}), an increase in the stability of the difference scheme can be achieved
twofold. In the first case, due to an increase in the energy $(Cy, y)$ of the operator
$C$ or by reducing the account
energy of the operator $D$. The first possibility of constructing stable difference schemes
is based on using additive regularization: increasing operator $C$ or/and decreasing operator
$D$ due to additional terms. The second possibility is related to the
multiplicative perturbation of the operators of the generating scheme.

With the multiplicative regularization of the operator $C$,
for example, we will replace $C \longmapsto C (I + \mu R)$ or $C \longmapsto (I + \mu R) C$,
where $R$ is a regularizing operator and $\mu$ is a regularization parameter.
With such a perturbation, we remain in the class of schemes with
self-adjoint operators if $R C = C R^*$.
An example of a more complex regularization is given by the transformation
\[
  C \longmapsto (I + \mu R^*) C (I + \mu R).
\]

The multiplicative regularization is carried out similarly due to the perturbation
operator $D$. Taking into account the inequality (\ref{2.9}), we can implement
transformation $D \longmapsto D (I + \mu R)^{-1} $ or $D \longmapsto (I + \mu R)^{-1} D$.
For the simplest two-level schemes, such a regularization can
consider as a new edition of the regularization of the operator $C$. 
To stay in the class of schemes with self-adjoint operators, it is enough to
choose $R = R(D)$. We have great opportunities for regularization
\[
  D \longmapsto (I + \mu  R^*)^{-1} D (I + \mu  R)^{-1}.
\]
In this case, the regularizing operator $R$ may not directly bind to
operator $D$.

Under perturbation of the operator $D$ from (\ref{3.2}), we arrive at the scheme
\begin{equation}\label{3.3}
  \frac{u^{n+1} - 2 u^n + u^{n-1}} {\tau^2} + \widetilde{Q} u^n = 0,
  \quad n = 1,\ldots,N-1 .
\end{equation}
For multiplicative regularization, we have, for example,
$\widetilde{Q} = \widetilde{R} Q$. In the simplest case $\widetilde{R} = (I + \mu Q)^{-1}$ from (\ref{3.3}) we obtain
a regularized scheme
\begin{equation}\label{3.4}
  \frac{u^{n+1} - 2 u^n + u^{n-1}} {\tau^2} + (I + \mu Q)^{-1} Q u^n = 0,
  \quad n = 1,\ldots,N-1 .
\end{equation}
The scheme (\ref{3.4}) we obtain with the additive regularization of the operator with the time derivative: $C \longmapsto C + \mu Q, \ C = I$. 

Checking the inequality (\ref{2.9}) gives that for
\[
 \mu = \sigma \tau^2, 
 \quad \sigma \geq \frac{1}{4} , 
\] 
the regularized scheme (\ref{3.4}), (\ref{2.6}) is stable.
This scheme is directly related to the conventional weighted scheme for equation (\ref{3.1}):
\[
  \frac{u^{n+1} - 2 u^n + u^{n-1}} {\tau^2} +  Q (\sigma u^{n+1} + (1 - 2\sigma) u^n + \sigma u^{n-1}) = 0 ,
\]
whose stability conditions are well known \cite{SamarskiiTheory,SamarskiiGulin1973}. 

In the case of an additive representation of the operator $Q$, stable splitting schemes can be constructed
based on the perturbation of the operator terms. Let in the equation (\ref{3.1})
\[
 Q = \sum_{\alpha =1}^{p} Q_\alpha ,
 \quad  Q_\alpha^* =  Q_\alpha \geq 0,
 \quad \alpha =1,\dots,p . 
\] 
Similarly (\ref{3.3}), (\ref{3.4}), we will use the scheme
\begin{equation}\label{3.5}
  \frac{u^{n+1} - 2 u^n + u^{n-1}} {\tau^2} + \sum_{\alpha =1}^{p}\widetilde{Q}_\alpha  u^n = 0,
  \quad n = 1,\ldots,N-1 ,  
\end{equation} 
wherein
\[
 \widetilde{Q}_\alpha = (I + \mu_\alpha  Q_\alpha )^{-1}  Q_\alpha ,
 \quad \alpha =1,\dots,p .  
\] 
In the simplest case of equal weights $\mu_\alpha, \ \alpha = 1, \dots, p,$
this additive scheme will be stable when
\[
 \mu_\alpha = \sigma_\alpha  \tau^2, 
 \quad \sigma_\alpha = \sigma \geq \frac{ p}{4} ,
 \quad \alpha =1,\dots,p .   
\] 
Thus, stability is ensured by increasing the weighting factors.

The implementation of the scheme (\ref{3.5}) can be carried out based on solving independent problems
\[
  \frac{u^{n+1}_\alpha  - 2 u^n + u^{n-1}} {\tau^2} + \widetilde{Q}_\alpha  u^n = 0,
 \quad \alpha =1,\dots,p ,
\] 
and determining the solution on a new layer in time according to the rule
\[
 u^{n+1} = \frac{1}{p} \sum_{\alpha =1}^{p} u^{n+1}_\alpha ,
  \quad n = 1,\ldots,N-1 . 
\] 
Such an organization of computations corresponds to the use of an additive-averaged scheme \cite{VabishchevichAdditive}. 

We separately note the possibilities of multiplicative regularization for problems with the product of operators.
Let in the equation (\ref{3.1}) $Q = A^*A > 0 $ and $A> 0$.
Standard multiplicative regularization when in (\ref{3.3})
\[
 \widetilde{Q} = (I + \mu  A^* A)^{-1}  A^* A,
\] 
maybe unacceptable due to the need to solve the problem
with the operator $(I + \mu_\alpha A^*A)$. Therefore, it makes sense to consider the option with
perturbation of each operator factor in $Q = A^*A$.
For example, put
\begin{equation}\label{3.6}
 \widetilde{Q} = (I + \mu  A^*)^{-1}  A^* A (I + \mu  A)^{-1} .
\end{equation} 

Under the conditions (\ref{3.6}), the inequality (\ref{2.9}) for $C = I $, $ D = \widetilde {Q}$ is satisfied if
\[
  (I + \mu  A^*)(I + \mu  A) \geq \frac{\tau^2}{4 }   A^* A .
\] 
Thus, it suffices to put
\begin{equation}\label{3.7}
 \mu = \sigma \tau, 
 \quad \sigma \geq \frac{1}{2} .
\end{equation} 
The main potential drawback of regularization (\ref{3.6}), (\ref{3.7}) is related to the fact that
\[
 \widetilde{Q} = Q + (A^* + A) \, \mathcal{O} (\tau) .
\] 
In case (\ref{3.4}) we have
\[
 \widetilde{Q} = Q + A^* A  \, \mathcal{O} (\tau^2) ,
\] 
that is, the perturbation is associated with the second-order in $\tau$.

\section{Regularized scheme} 

Now we can construct an unconditionally stable scheme based on the multiplicative regularization
for our problem (\ref{2.1})--(\ref{2.3}).
It is associated with the perturbation of the operators $A^*A$ and $B$ and has the form
\begin{equation}\label{4.1}
  \frac{u^{n+1} - 2 u^n + u^{n-1}} {\tau^2} + \widetilde{A^* A}  u^n + \widetilde{B}  u^n = 0 ,
  \quad n = 1,\ldots,N-1 .
\end{equation}
This scheme is written in the form (\ref{2.7}) with
\[
 C = I,
 \quad D =  \widetilde{A^* A} + \widetilde{B} .
\] 
According to Lemma \ref{l-1}, stability will be ensured, in particular, for
\begin{equation}\label{4.2}
 I \geq \frac{\tau^2}{2} \widetilde{A^* A},
 \quad   I \geq \frac{\tau^2}{2} \widetilde{B} .
\end{equation} 

For the regularizing operator $\widetilde{B}$ put
\begin{equation}\label{4.3}
 \widetilde{B} = (I + \sigma_B \tau^2 B)^{-1} B ,
 \quad  \sigma_B > 0 .
\end{equation} 
We will perturb the operator $A^*A$ according to (\ref{3.6}):
\begin{equation}\label{4.4}
 \widetilde{A^* A} = (I + \sigma_A \tau A^*)^{-1} A^* A (I + \sigma_A \tau A)^{-1} ,
 \quad  \sigma_A > 0 . 
\end{equation} 
For such $\widetilde{A^*A} $ and $\widetilde{B}$ the inequalities (\ref{4.2}) will hold for
the following restrictions on weight parameters:
\begin{equation}\label{4.5}
 \sigma_A^2 \geq \frac{1}{2} ,
 \quad  \sigma_B \geq \frac{1}{2} .
\end{equation} 
The result of our consideration is the following statement.

\begin{thm}\label{t-3}
The regularized scheme (\ref{2.6}), (\ref{4.1}), (\ref{4.3}), (\ref{4.4}) is unconditionally stable for
constraints (\ref{4.5}).
\end{thm}

Similarly (\ref{3.5}), schemes with additional splitting are constructed.
The simplest variant is associated with splitting operator $B$. Let in (\ref{2.3})
\[
 B = \sum_{\alpha =1}^{p} B_\alpha ,
 \quad  B_\alpha^* =  B_\alpha \geq 0,
 \quad \alpha =1,\dots,p . 
\] 
In the scheme (\ref{4.1}), (\ref{4.4}) we define
\begin{equation}\label{4.6}
 \widetilde{B} = \sum_{\alpha =1}^{p}(I + \sigma_B \tau^2 B_\alpha )^{-1} B_\alpha .
\end{equation}  
The stability of the scheme (\ref{2.6}), (\ref{4.1}), (\ref{4.4}), (\ref{4.6}) takes place, for example, for
\[
 \sigma_A^2 \geq \frac{1}{2} ,
 \quad  \sigma_B \geq \frac{p}{2} .
\] 

Similarly, we consider the case of splitting the operator $A^*A$, when
\[
 A^* A = \sum_{\alpha =1}^{p} A^*_\alpha A_\alpha ,
 \quad \widetilde{A^* A} = \sum_{\alpha =1}^{p}  (I + \sigma_A \tau A_\alpha^*)^{-1} A_\alpha^* A_\alpha(I + \sigma_A \tau A_\alpha)^{-1} .
\] 
In this case, the scheme (\ref{2.6}), (\ref{4.1}), (\ref{4.3}) is unconditionally stable for
\[
 \sigma_A^2 \geq \frac{p}{2} ,
 \quad  \sigma_B \geq \frac{1}{2} .
\] 

The variant of splitting the operator's $A$ and $A^*$ deserves special attention, when
\[
 A = \sum_{\alpha =1}^{p} A_\alpha ,
 \quad  A^* = \sum_{\alpha =1}^{p} A^*_\alpha,
 \quad \alpha =1,\dots,p .  
\]  
For $\widetilde {A^* A} $ put
\begin{equation}\label{4.7}
 \widetilde{A^* A} = \sum_{\alpha =1}^{p} \widetilde{A}^*_\alpha \sum_{\alpha =1}^{p} \widetilde{A}_\alpha ,
 \quad  \widetilde{A}_\alpha = (I + \sigma_A \tau A_\alpha)^{-1} A_\alpha,
 \quad  \alpha =1,\dots,p .   
\end{equation} 
In this case, we have
\[
 (\widetilde{A^* A} u, u) = \left ( \Big (\sum_{\alpha =1}^{p} \widetilde{A}_\alpha u \Big )^2, 1 \right ) 
 \leq p  \sum_{\alpha =1}^{p} \left (\Big ( \widetilde{A}_\alpha u \Big )^2, 1 \right ) =
 p  \sum_{\alpha =1}^{p} \Big (\widetilde{A}^*_\alpha \widetilde{A}_\alpha u, u \Big ) .
\] 
With this in mind, inequalities (\ref{4.2}) will be satisfied with
\[
 \sigma_A^2 \geq \frac{p^2}{2} ,
 \quad  \sigma_B \geq \frac{1}{2} ,
\] 
which ensures the stability of the scheme   (\ref{2.6}), (\ref{4.1}),  (\ref{4.3}),  (\ref{4.7}).

\section{Numerical experiments} 

The possibilities of using the constructed splitting schemes will be illustrated by the results of the numerical solution
model two-dimensional problem. We will assume that the computational domain is a rectangle
\[
 \Omega = \{ \bm x  \ | \ \bm x = (x_1,x_2), \ 0 < x_\alpha  < l_\alpha , \ \alpha  = 1,2 \} ,
\]
with boundary $\partial \Omega$.
We need to find a solution $v(\bm x, t)$ of the equation 
\begin{equation}\label{5.1}
 \frac{\partial^2 v}{\partial t^2} +  \triangle^2 v + \gamma_1 v - \gamma_2 \triangle v = 0,
 \quad \bm x \in \Omega ,
 \quad 0 < t \leq T ,  
\end{equation} 
where $\gamma_1 = \const > 0$,  $\gamma_2 = \const > 0$, and $\triangle = \div \grad$ is the Laplace operator.
Equation (\ref{5.1}) is supplemented with the following boundary and initial conditions:
\begin{equation}\label{5.2}
 v(\bm x,t) = 0,
 \quad  \triangle v(\bm x,t) = 0,
 \quad \bm x \in \partial \Omega , 
\end{equation} 
\begin{equation}\label{5.3}
 v(\bm x,0) = v^0(\bm x), 
 \quad \frac{\partial v}{\partial t} (\bm x,0) = 0, 
 \quad \bm x \in \Omega .
\end{equation} 

Boundary value problem (\ref{5.1})--(\ref{5.3}) describes (see details, for example, in the \cite{selvadurai1979elastic})
displacement of the plate on the elastic base. In this case, $v(\bm x, t) $ is the normal displacement of platinum,
$v^0 (\bm x)$ defines the offset at the start time. The boundary conditions (\ref{5.2}) correspond to
hinge fastening. In the framework of two-dimensional elastic models, the parameter $\gamma_1$ is associated with
the elastic foundation reaction modulus (Winkler model), and the $\gamma_2 $ parameter --- with the tension action of a thin elastic membrane
in the Filonenko-Borodich model and with the shear action among the spring elements in the Pasternak model.

On the set of sufficiently smooth functions $w(\bm x) = 0, \ \bm x \in \partial \Omega$, we define the operator
\[
 \mathcal{A} w = - \triangle w, 
 \quad \bm x \in \Omega .
\] 
Let us write the problem (\ref{5.1})--(\ref{5.3}) in the form of the Cauchy problem for a second-order evolution equation.
The solution $v(t) = v(\cdot, t)$ is determined from the equation
\begin{equation}\label{5.4}
 \frac{d^2 v}{d t^2} + \mathcal{A}^2 v + \gamma_1 v + \gamma_2 \mathcal{A} v = 0,
 \quad 0 < t \leq T .
\end{equation} 
Taking into account (\ref{5.3}), it is supplemented with the initial conditions
\begin{equation}\label{5.5}
 v(0) = v^0,
 \quad \frac{d v}{d t} (0) = 0 .
\end{equation} 

To numerically solve the problem (\ref{5.1})--(\ref{5.3}), we will use the standard
difference approximations in space \cite{SamarskiiTheory}.
We will introduce in the region $\Omega$ a uniform rectangular grid
\[
\overline{\omega}  = \{ \bm{x} \ | \ \bm{x} =\left(x_1, x_2\right), \quad x_\alpha  =
i_\alpha  h_\alpha , \quad i_\alpha  = 0,1,...,N_\alpha ,
\quad N_\alpha  h_\alpha  = l_\alpha , \ \alpha  = 1,2 \} ,
\]
where $\overline{\omega} = \omega \cup \partial \omega$, and
$\omega$ is the set of internal ones, and $\partial \omega$ is the set of boundary mesh nodes.
For grid functions $w(\bm x)$ such that $w(\bm x) = 0, \ \bm x \notin \omega$, we define the Hilbert space
$H = L_2 (\omega)$, in which the dot product and norm are
\[
(w, u) = \sum_{\bm x \in  \omega} w(\bm{x}) u(\bm{x}) h_1 h_2,  \quad 
\| w \| =  (w,w)^{1/2}.
\]

For $u(\bm x) = 0, \ \bm x \notin \omega$, we define the grid Laplace operator $A$ on the usual five-point pattern:
\[
  \begin{split}
  A u = & -
  \frac{1}{h_1^2} (u(x_1+h_1,x_2) - 2u(\bm{x}) + u(x_1-h_1,x_2)) \\ 
  & - \frac{1}{h_2^2} (u(x_1,x_2+h_2) - 2 u(\bm{x}) + u(x_1,x_2-h_2)), 
  \quad \bm{x} \in \omega . 
 \end{split} 
\] 
For this grid operator (see, for example, \cite{SamarskiiTheory}) we have
\begin{equation}\label{5.6}
 A = A^* \geq \delta I,
 \quad \delta > 0 . 
\end{equation} 
On sufficiently smooth functions, the operator $A$ approximates the differential operator
$\mathcal{A}$ with an error
$\mathcal{O} \left(|h|^2\right)$, $|h|^2 = h_1^2+h_2^2$. 

The finite-difference approximation in the space of the problem (\ref{5.4}), (\ref{5.5}) leads us to
equation (\ref{2.1}),
which is complemented by the initial conditions
\begin{equation}\label{5.7}
 w(0)= w^0 ,
 \quad \frac{d w}{d t} (0) = 0 ,
\end{equation}
when $w^0 = v^0(\bm x) = 0, \ \bm x \in \omega$.
For the operator $B$ we have
\begin{equation}\label{5.8}
 B = \gamma_1 I + \gamma_2 A . 
\end{equation} 

We carry out numerical experiments based on the exact solution of the problem (\ref{2.1}), (\ref{5.6})--(\ref{5.8}).
Consider the grid spectral problem
\[
 A \psi =  \lambda  \psi . 
\] 
For eigenfunctions and eigenvalues we have (see, for example, \cite{SamarskiiNikolaev1978}):
\[
\begin{split}
 \psi_k (\bm x) & = \prod_{\beta  =1}^{2} \sqrt{\frac{2}{l_\beta } } \sin (k_\beta  \pi x_\beta ),
 \quad \bm x \in \omega , \\ 
 \lambda_k & = \sum_{\beta  =1}^{2}  \frac{4}{h_\beta^2} \sin^2 \frac{k_\beta  \pi}{2 N_\beta } ,
 \quad k_\alpha  = 1,2,...,N_\alpha -1, \quad  \alpha  = 1, 2 .
\end{split} 
\] 
Because of this
\[
 \delta = \lambda _1 = \sum_{\beta  =1}^{2}  \frac{4}{h_\beta ^2} \sin^2 \frac{\pi}{2 N_\beta } 
 <  8 \Big (\frac{1}{l_1^2} + \frac{1}{l_2^2} \Big ) .
\] 
Eigenfunctions $\psi_k, \ \|\psi_k \| = 1, \,$ form a basis in $H$.
Therefore, for any grid function $u \in H$, the representation takes place
\[
 u = \sum_{k= 1}^{K}(u, \psi_k) \psi_k . 
\]
To solve the problem (\ref{2.1}), (\ref{5.6})--(\ref{5.8}) we get
\begin{equation}\label{5.9}
 w(\bm x, t) = \sum_{k= 1}^{K}(w^0, \psi_k) \cos(r_k^{1/2} t) \psi_k (\bm x),
 \quad r_k =  \gamma_1 + \gamma_2 \lambda_k + \lambda_k^2 .
\end{equation} 

\begin{figure}
\centering
\begin{minipage}{0.49\linewidth}
\centering
\includegraphics[width=\linewidth]{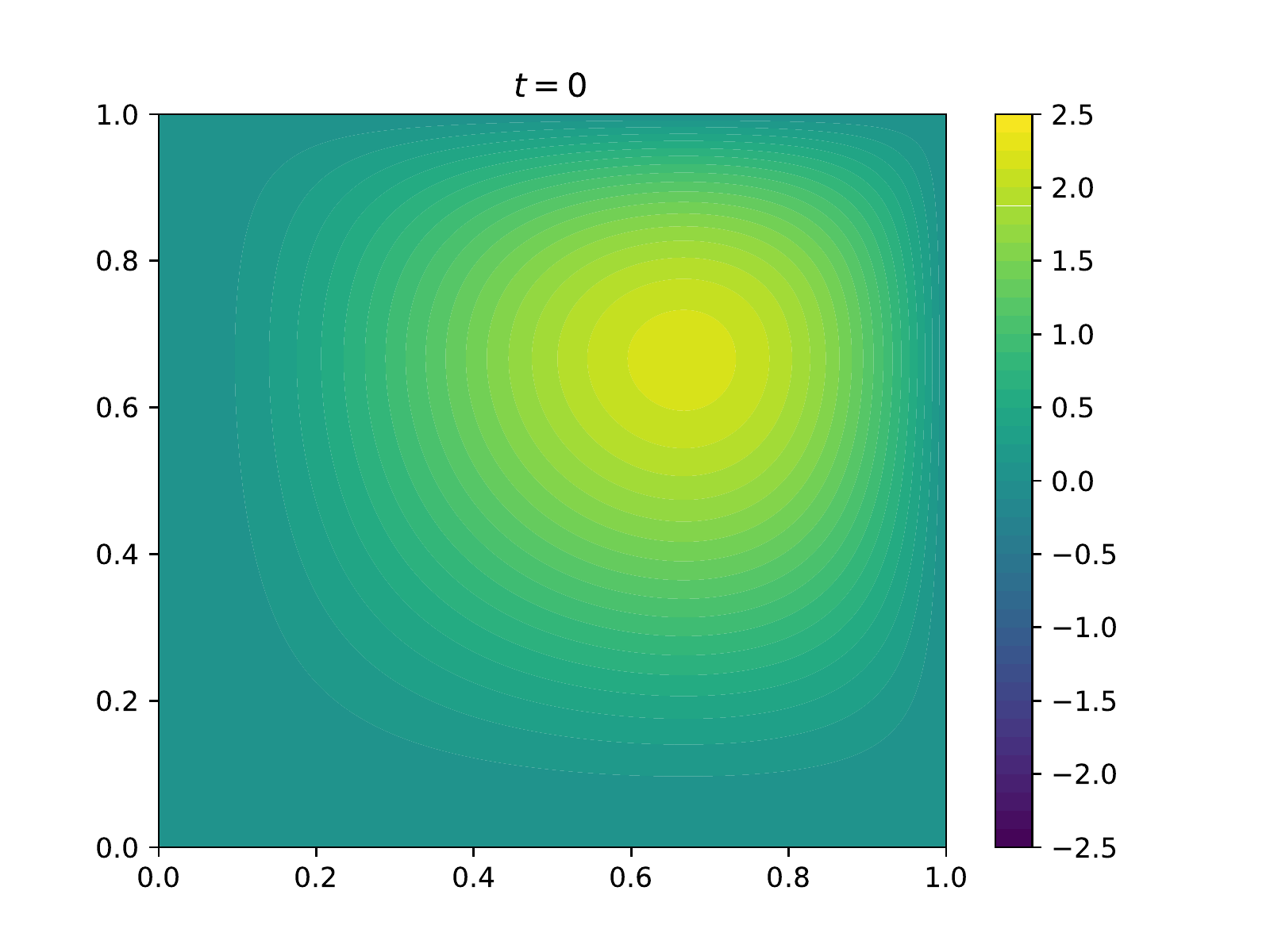}\\
\includegraphics[width=\linewidth]{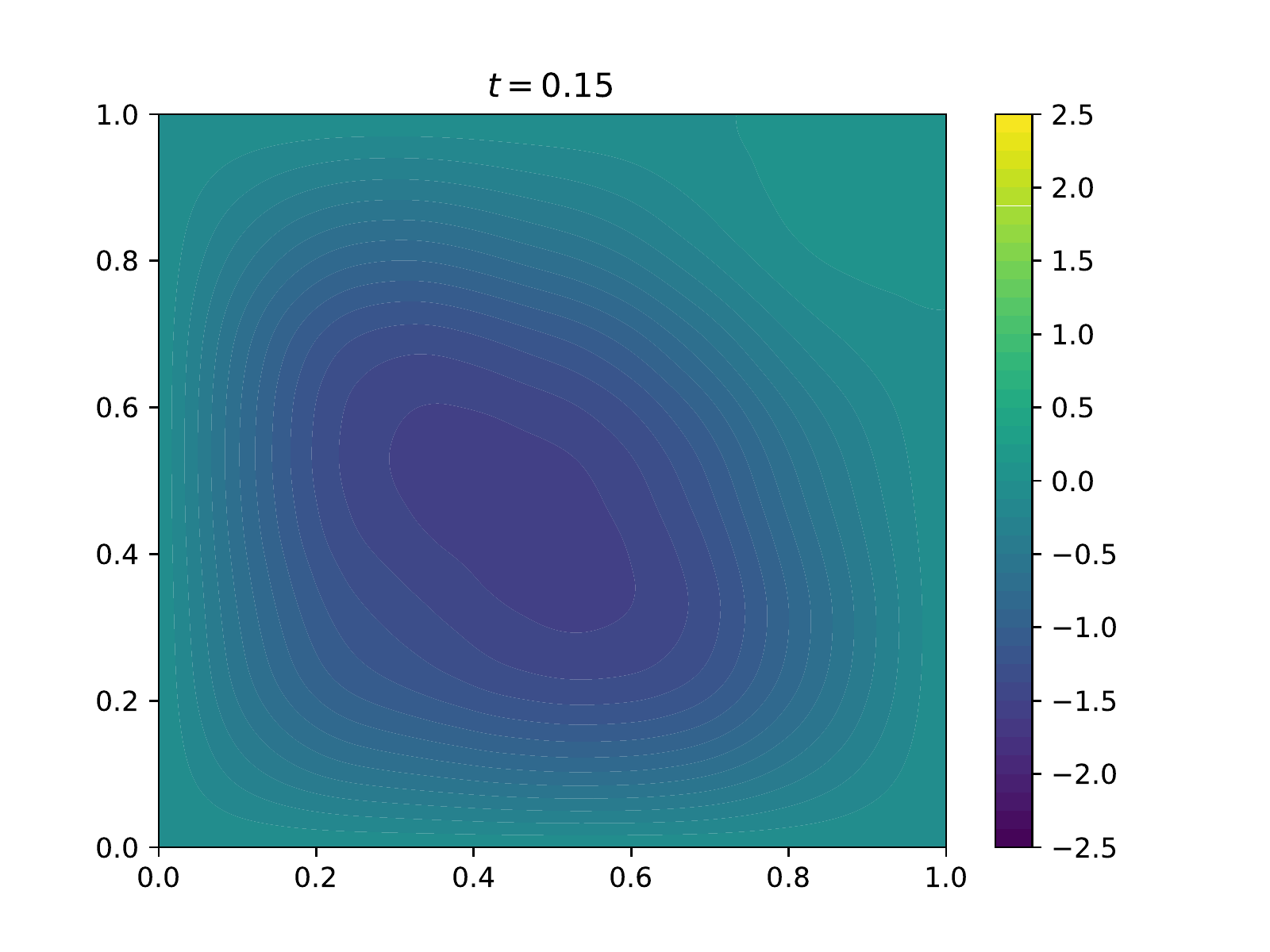}\\
\end{minipage}
\begin{minipage}{0.49\linewidth}
\centering
\includegraphics[width=\linewidth]{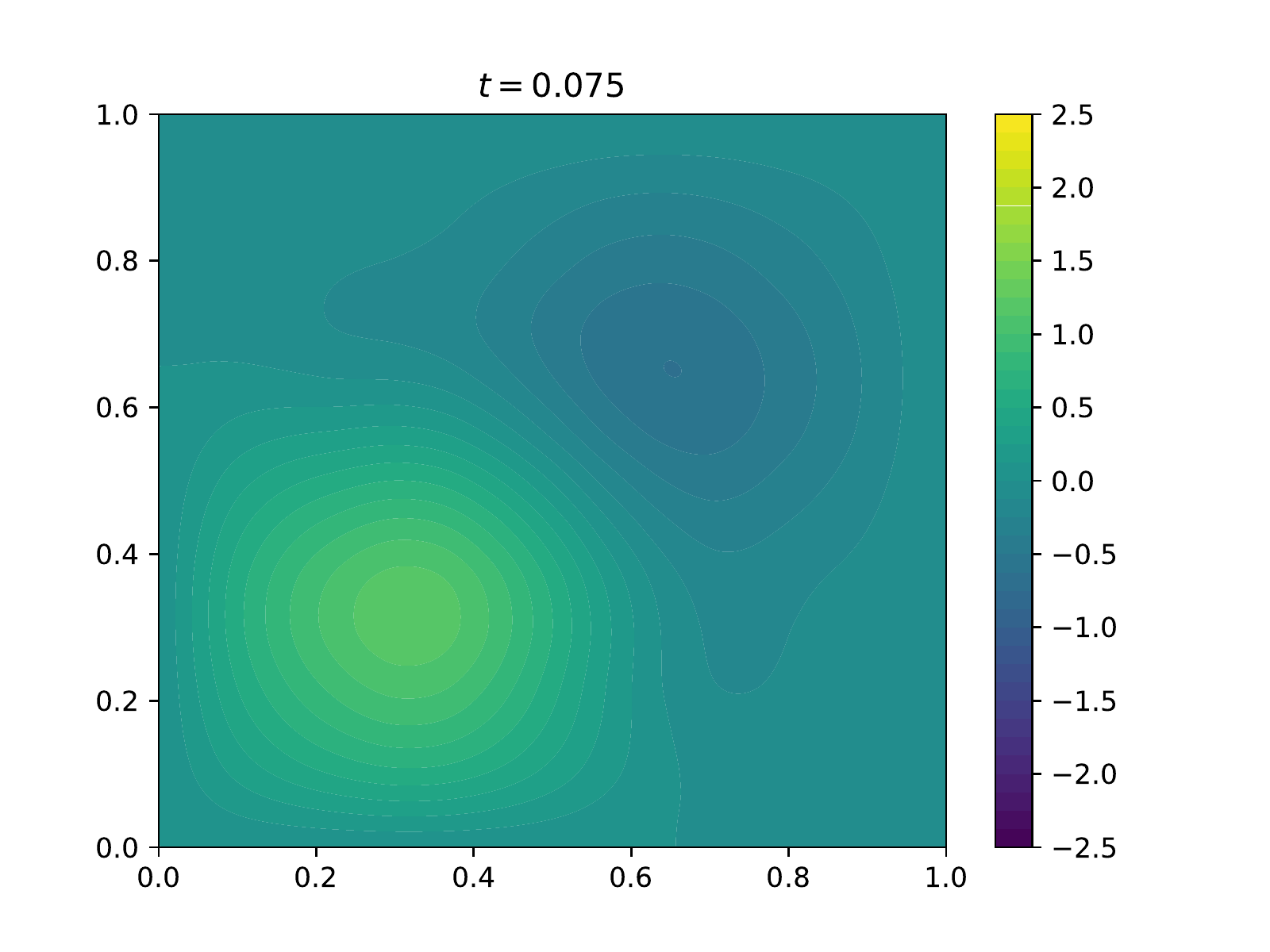}\\
\includegraphics[width=\linewidth]{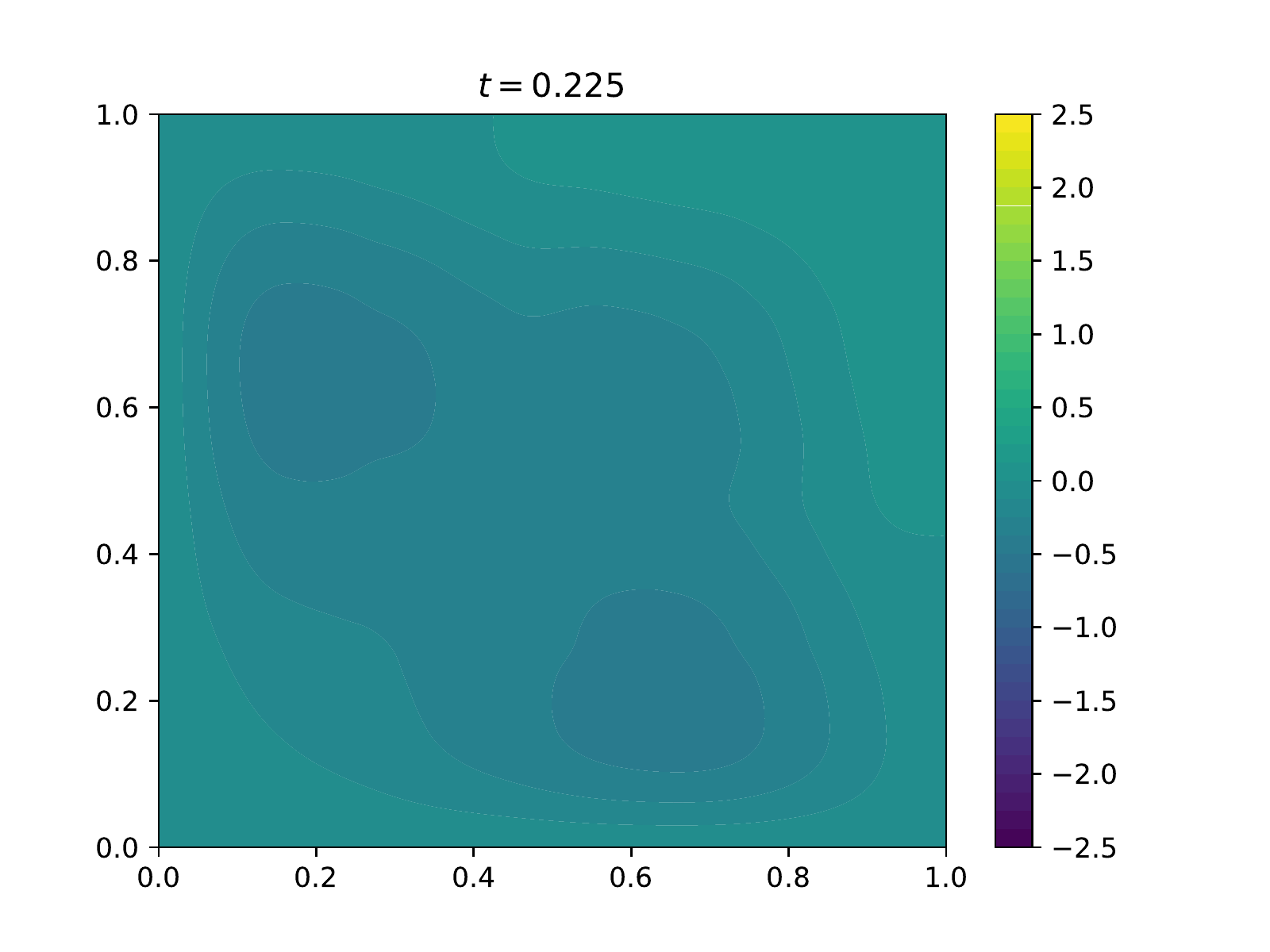}\\
\end{minipage}
\caption{Exact solution of the problem at separate points in time.}
\label{f-1}
\end{figure}

\begin{figure}
\centering
\includegraphics[width=0.75\linewidth]{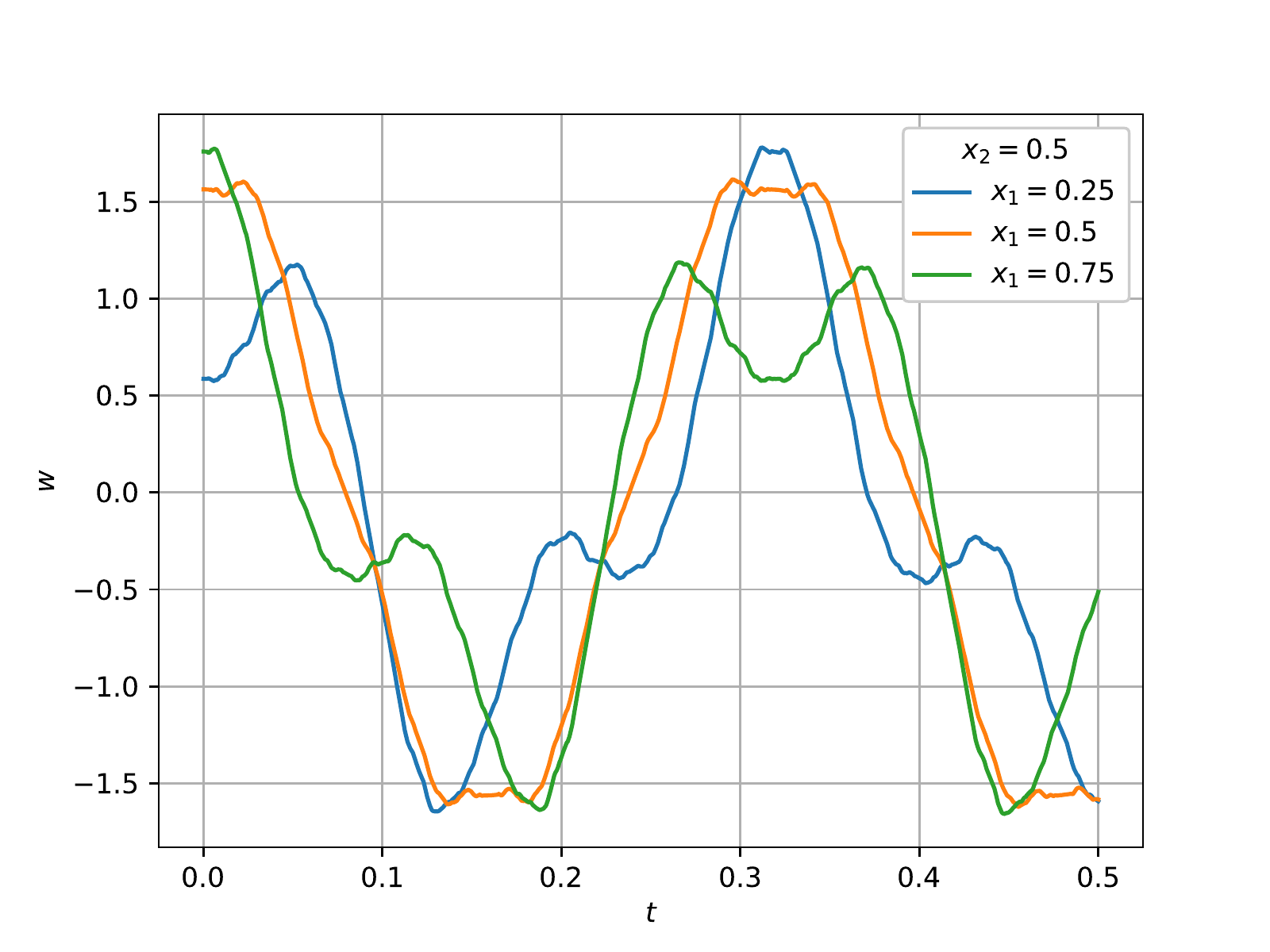}
\centering
\caption{Deflection of the plate at individual points.}
\label{f-2}
\end{figure}

The calculation results presented below were obtained for the problem with
\[
 l_1 = l_2 = 1, 
 \quad N_1 = N_2 = 256,
 \quad \gamma_1 = 1, 
 \quad \gamma_2 = 0.05 .   
\] 
Comparison of approximations in time is carried out on the problem with the initial condition
\[
 w^0(\bm x) = x_1^2(1-x_1) x_2^2(1-x_2).
\] 
The solution to the test problem at separate points in time is shown in Fig.\ref{f-1}.
The plate deflection dynamics at characteristic points are shown in Fig.\ref{f-2}.
Complex vibrations of the plate are observed.

\begin{figure}
\centering
\begin{minipage}{0.49\linewidth}
\centering
\includegraphics[width=\linewidth]{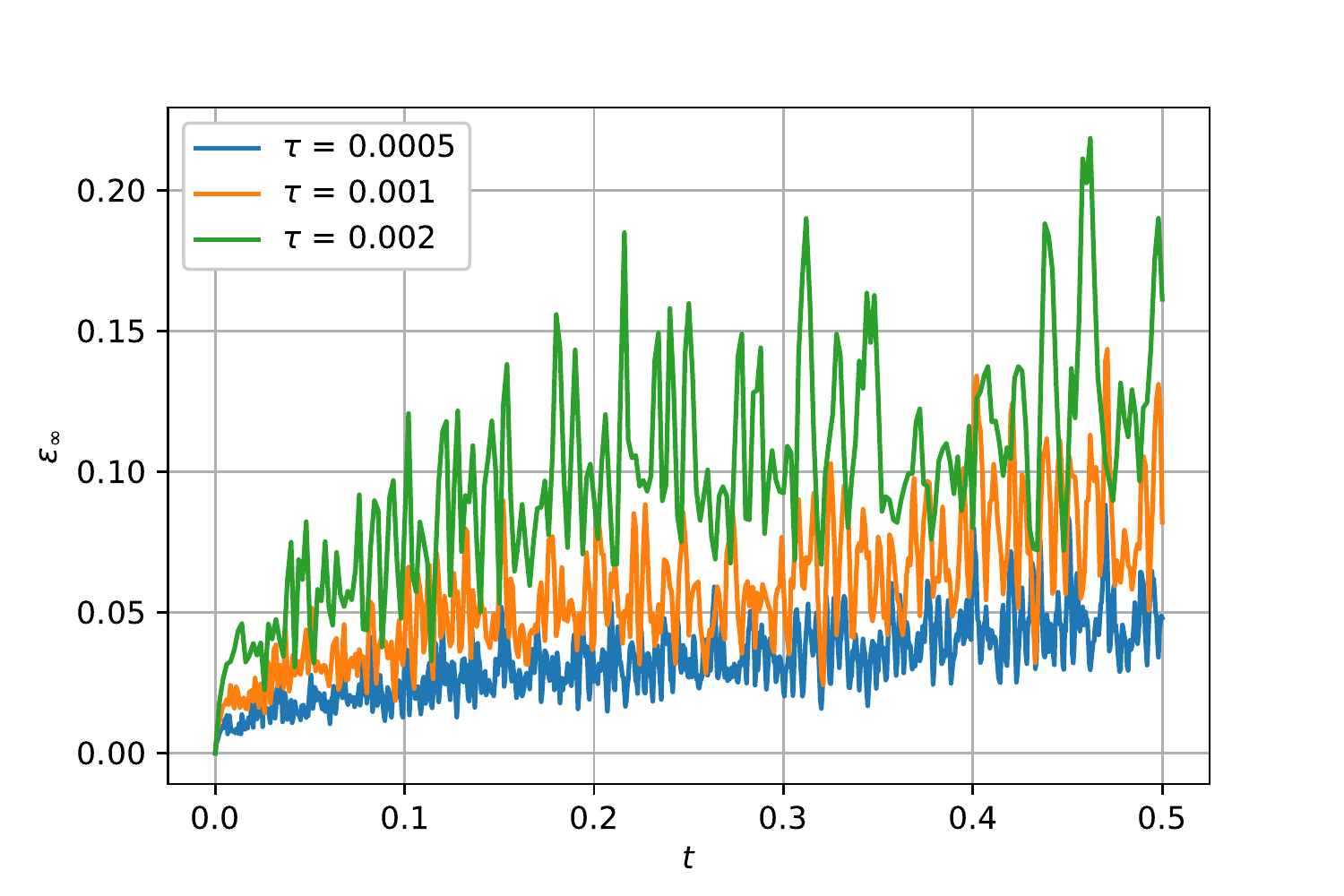}\\
\end{minipage}
\begin{minipage}{0.49\linewidth}
\centering
\includegraphics[width=\linewidth]{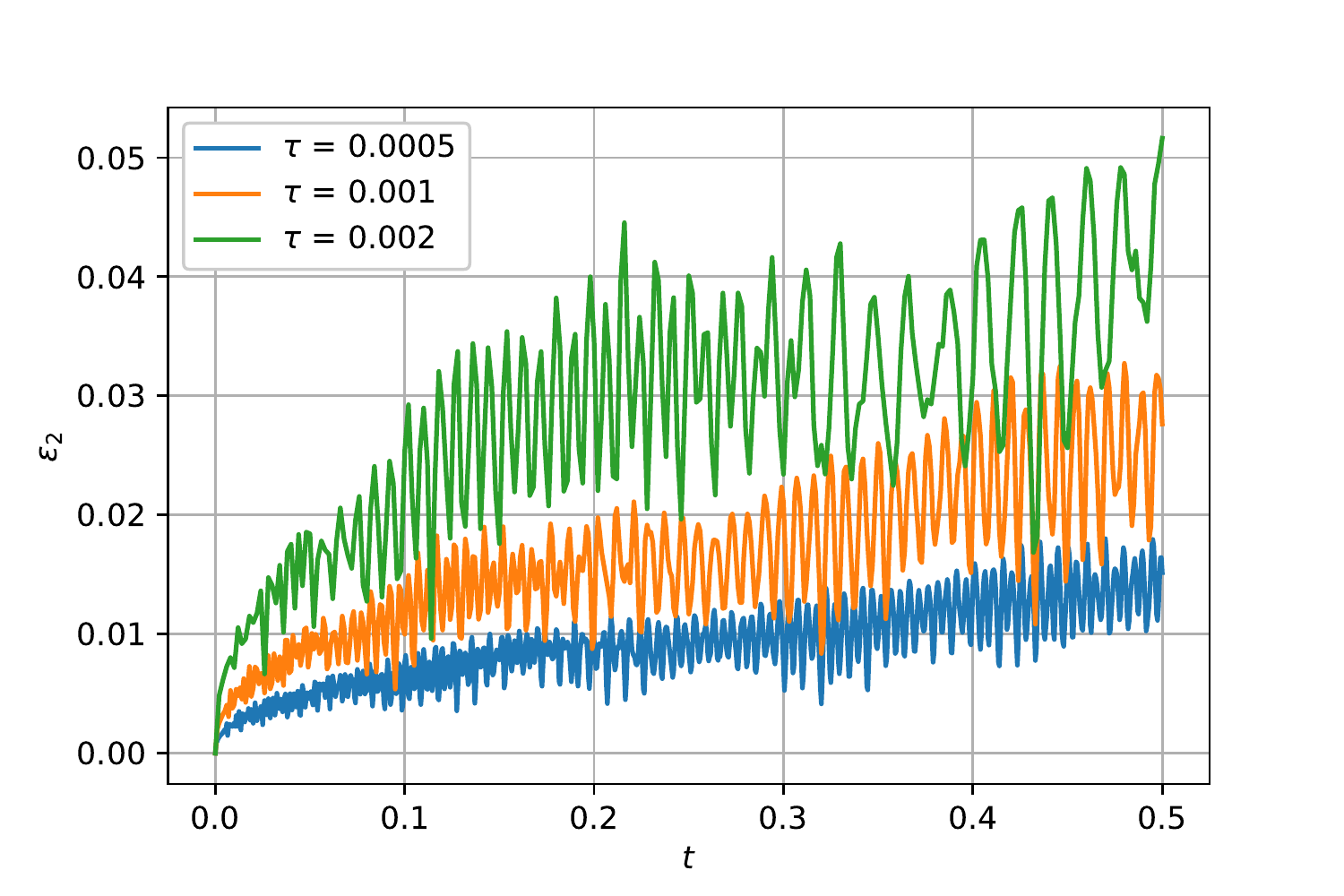}\\
\end{minipage}
\caption{The solution error for the weighted scheme at $\sigma = 0.25$.}
\label{f-3}
\end{figure}

\begin{figure}
\centering
\begin{minipage}{0.49\linewidth}
\centering
\includegraphics[width=\linewidth]{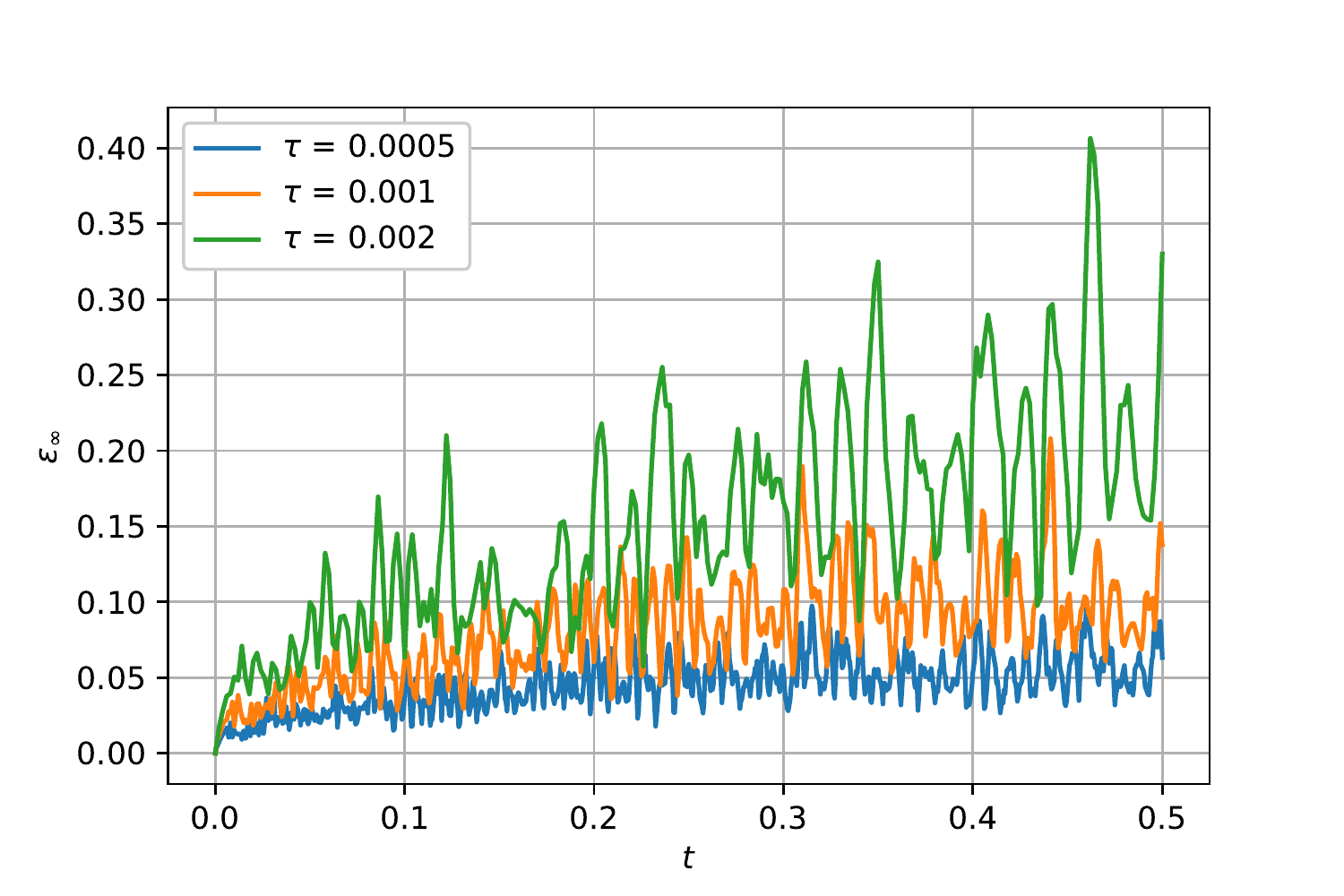}\\
\end{minipage}
\begin{minipage}{0.49\linewidth}
\centering
\includegraphics[width=\linewidth]{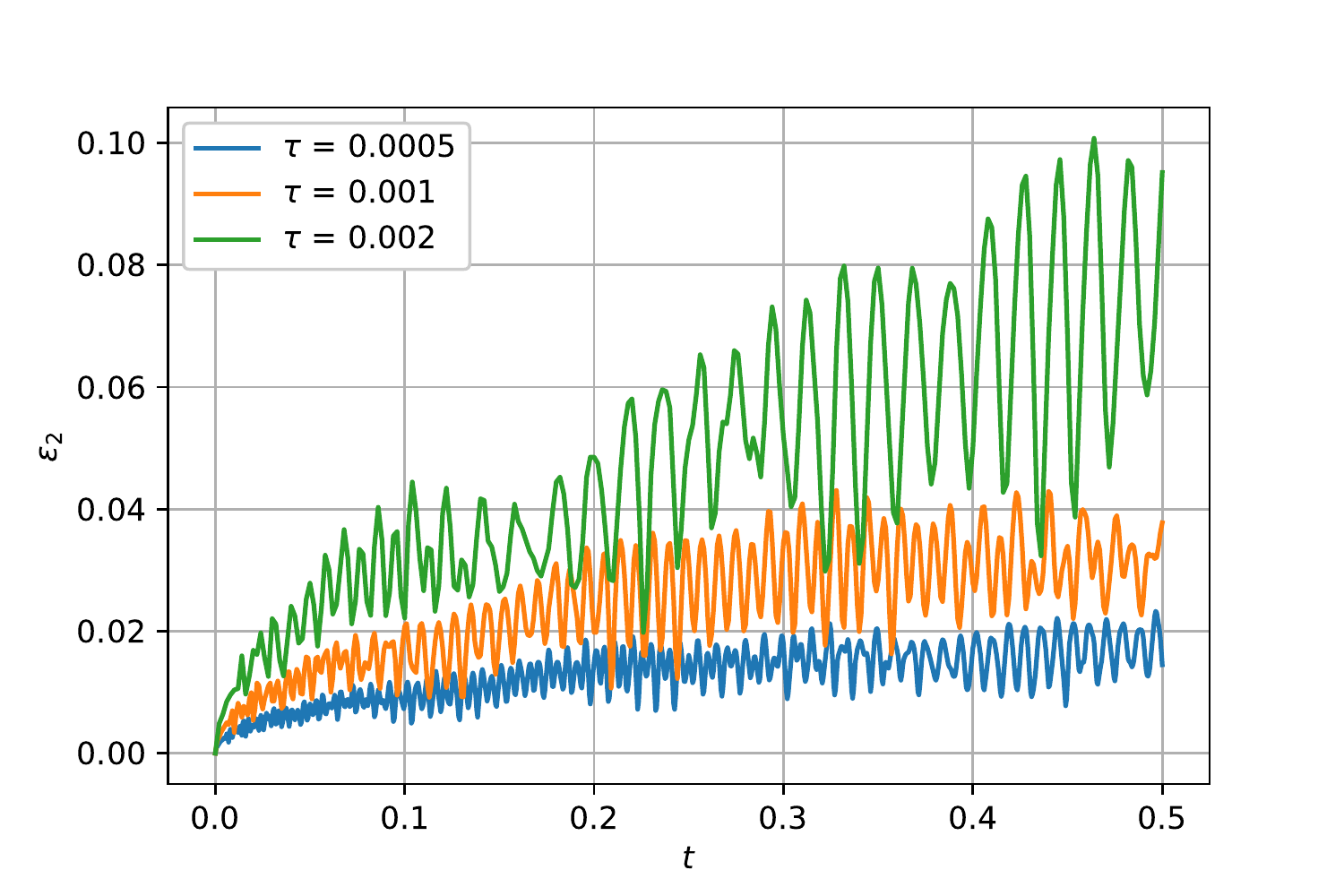}\\
\end{minipage}
\caption{The solution error for the weighted scheme at $\sigma = 0.5$.}
\label{f-4}
\end{figure}

It is natural to compare the accuracy of the constructed splitting scheme with the accuracy
usually weighted scheme (\ref{2.5}), (\ref{2.6}). At separate times $t = t^n$
we define the error norms of the solution in $C(\omega)$ and $L_2(\omega)$:
\[
 \varepsilon_\infty (t^n) = \max_{\bm x \in \omega} |u^n(\bm x) - w(\bm x, t^n) |,
 \quad  \varepsilon_2 (t^n) = \|u^n(\bm x) - w(\bm x, t^n) \|,
 \quad  n = 0, \ldots, N .
\]   
For the initial deflection of the plate, we have
\[
 \max_{\bm x \in \omega} |w(\bm x, 0) | \approx 2.195,
 \quad  \|w(\bm x, 0) \| \approx 0.9524 .
\] 

The error in the approximate solution of the problem (\ref{5.1})--(\ref{5.3}) is shown in Fig.\ref{f-3} when using the scheme (\ref{2.5}), (\ref{2.6}) with $\sigma = 0.25$.
For the considered initial data and time steps, the theoretical asymptotic dependence of the time step accuracy (second-order)
not visible. With increasing weight $\sigma$ the errors grow --- see Fig.\ref{f-4}.

\begin{figure}
\centering
\begin{minipage}{0.49\linewidth}
\centering
\includegraphics[width=\linewidth]{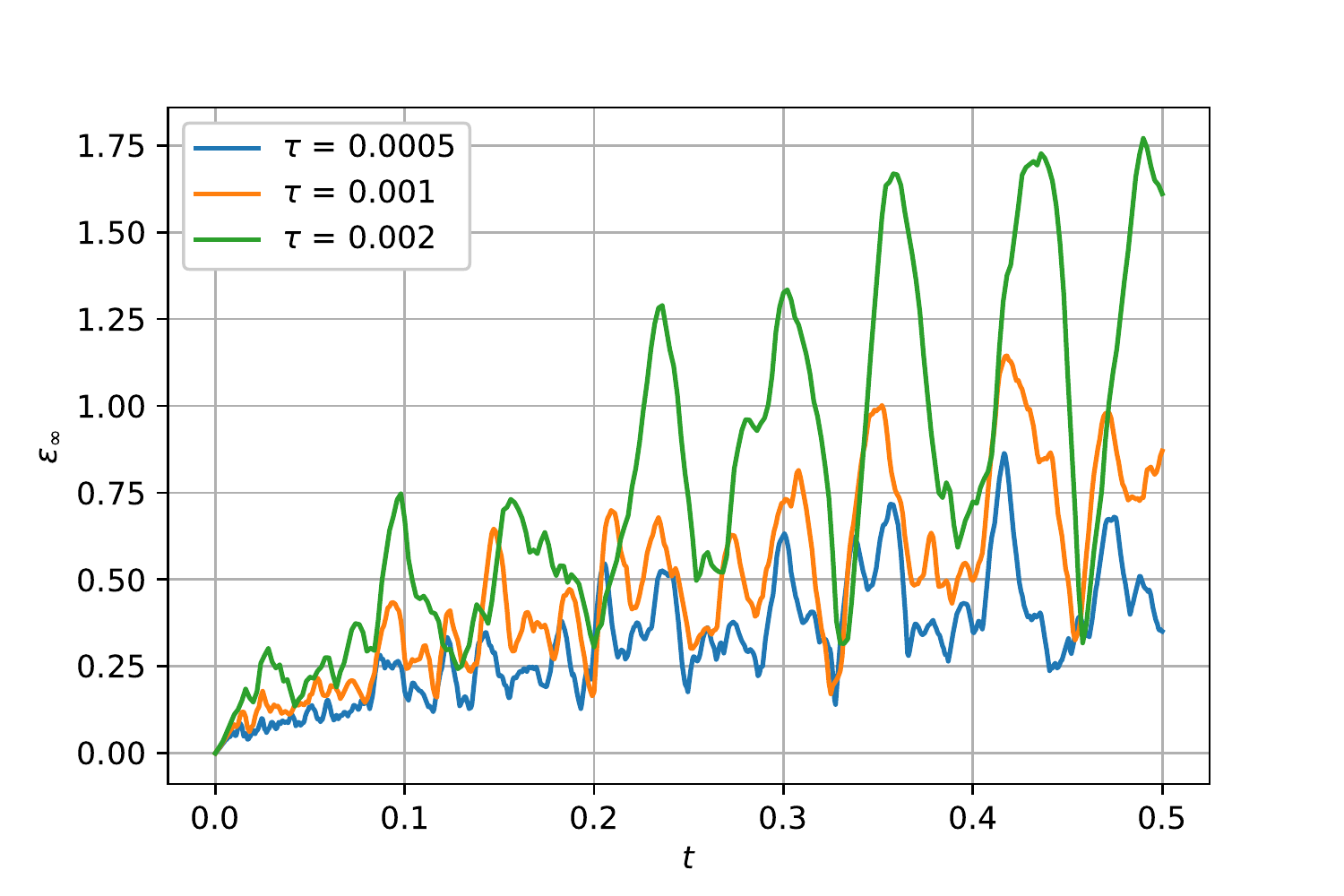}\\
\end{minipage}
\begin{minipage}{0.49\linewidth}
\centering
\includegraphics[width=\linewidth]{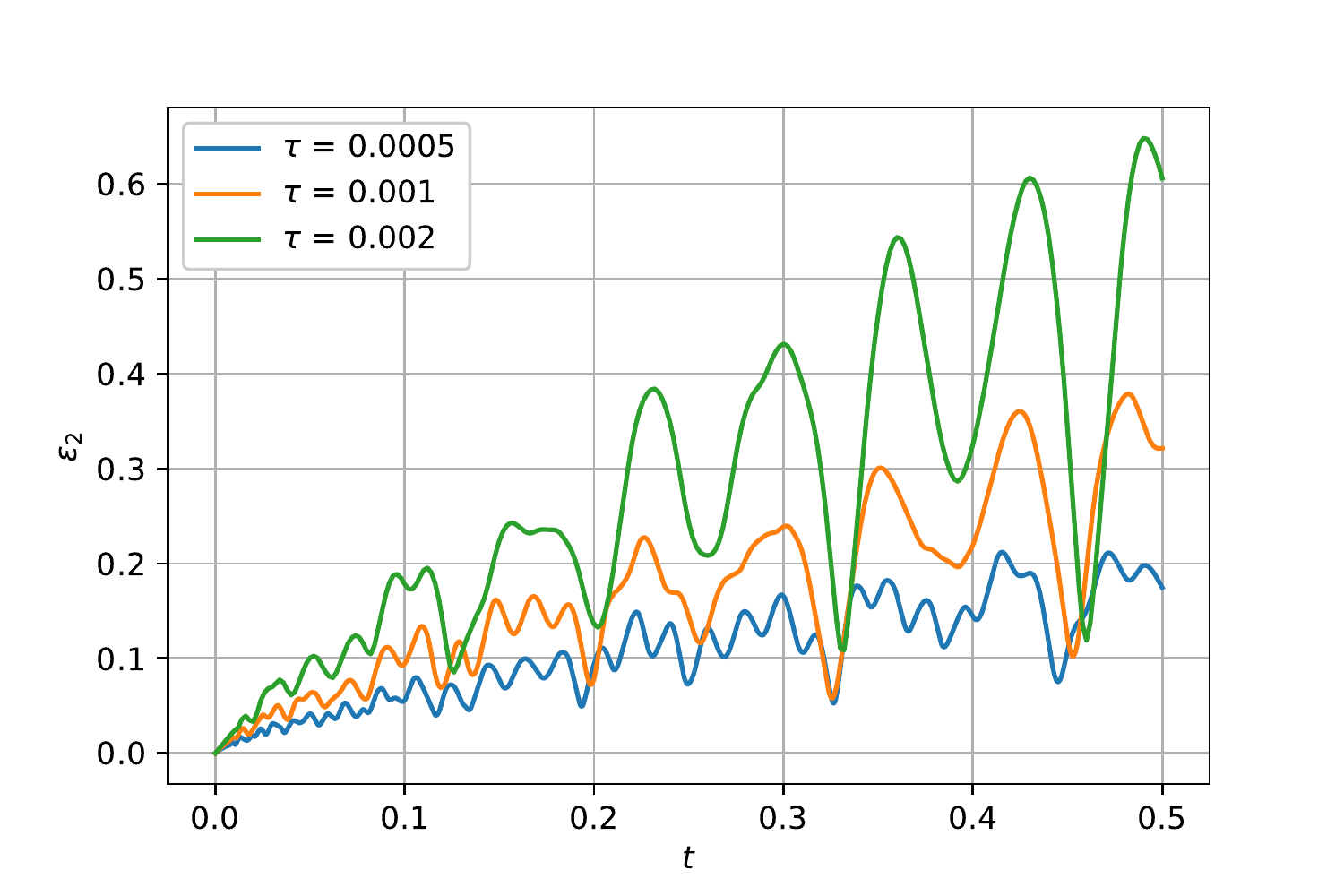}\\
\end{minipage}
\caption{The solution error for the splitting scheme with $\sigma_A^2 = 0.5, \sigma_B = 0.5$.}
\label{f-5}
\end{figure}

When using the splitting scheme (\ref{2.6}), (\ref{4.1}), (\ref{4.3}), (\ref{4.4})
subject to the constraints (\ref{4.5}), we set
\[
 \sigma_A^2 = \frac{1}{2} ,
 \quad  \sigma_B = \frac{1}{2} .
\] 
The dependence of the error of the approximate solution on time for this case is shown in Fig.\ref{f-5}.
Of course, as the time step decreases, the accuracy increases, but, as you would expect,
errors in comparison with the usual scheme with weights (see Fig.\ref{f-3}, \ref{f-4}) are much larger. 

\section{Conclusions} 

Applied models of the theory of plates lead to the necessity of solving the initial boundary value problems
for partial differential equations that include fourth-order elliptic operators.
The paper discusses the problem of reducing the computational complexity of the implementation
of unconditionally stable implicit schemes for these problems due to the use of special approximations in time.
In that work:
\begin{enumerate}
 \item A class of Cauchy problems for a second-order evolution equation is distinguished,
in which the problem operator is the sum of two self-adjoint operators. Wherein
one of the operators is represented as
the product of the operator $A$ by its conjugate $A^*$.
 \item Conditions for the absolute stability of three-level schemes are given.
with weights using general results
the theory of stability (correctness) of operator-difference schemes.
 \item Splitting schemes are constructed and investigated for which
the transition to a new level in time is associated with a separate solution of problems for
operators $A$ and $A^*$, not their products.
 \item The efficiency of the proposed splitting schemes is demonstrated by
the results of calculations of the dynamics of a thin square plate on an elastic foundation.
\end{enumerate}

\section*{Acknowledgements}

The publication has been prepared with support of the mega-grant of
the Russian Federation Government 14.Y26.31.0013.


\begin{thebibliography}{10}
\expandafter\ifx\csname url\endcsname\relax
  \def\url#1{\texttt{#1}}\fi
\expandafter\ifx\csname urlprefix\endcsname\relax\def\urlprefix{URL }\fi
\expandafter\ifx\csname href\endcsname\relax
  \def\href#1#2{#2} \def\path#1{#1}\fi

\bibitem{fung2017classical}
Y.-C. Fung, P.~Tong, X.~Chen, Classical and Computational Solid Mechanics,
  World Scientific, New Jersey, 2017.

\bibitem{SamarskiiTheory}
A.~A. Samarskii, The theory of difference schemes, Marcel Dekker, New York,
  2001.

\bibitem{LeVeque2007}
R.~J. LeVeque, Finite Difference Methods for Ordinary and Partial Differential
  Equations. Steady-State and Time-Dependent Problems, Society for Industrial
  Mathematics, 2007.

\bibitem{SamarskiiGulin1973}
A.~A. Samarskii, A.~V. Gulin, Stability of Difference Schemes, Nauka, Moscow,
  1973, in Russian.

\bibitem{SamarskiiMatusVabischevich2002}
A.~A. Samarskii, P.~P. Matus, P.~N. Vabishchevich, Difference Schemes with
  Operator Factors, Kluwer Academic Pub, 2002.

\bibitem{ascher1995implicit}
U.~M. Ascher, S.~J. Ruuth, B.~T.~R. Wetton, Implicit-explicit methods for
  time-dependent partial differential equations, SIAM Journal on Numerical
  Analysis 32~(3) (1995) 797--823.

\bibitem{HundsdorferVerwer2003}
W.~H. Hundsdorfer, J.~G. Verwer, Numerical solution of time-dependent
  advection-diffusion-reaction equations, Springer Verlag, 2003.

\bibitem{Marchuk1990}
G.~I. Marchuk, Splitting and alternating direction methods, in: P.~G. Ciarlet,
  J.-L. Lions (Eds.), Handbook of Numerical Analysis, Vol. I, North-Holland,
  1990, pp. 197--462.

\bibitem{VabishchevichAdditive}
P.~N. Vabishchevich, Additive Operator-Difference Schemes: Splitting Schemes,
  Walter de Gruyter GmbH, Berlin, Boston, 2013.

\bibitem{SamVabGul}
A.~A. Samarskii, P.~N. Vabishchevich, A.~V. Gulin, Stability of
  operator-difference schemes, Differential Equations 35~(2) (1999) 151--186.

\bibitem{VabFact}
P.~N. Vabishchevich, Operator-difference scheme with a factorized operator, in:
  Large-Scale Scientific Computing 10th International Conference, Sozopol,
  Bulgaria, June 8–12, 2015, Springer, 2015, pp. 72--79.

\bibitem{SamarskiiReg}
A.~A. Samarskii, Regularization of difference schemes, USSR Computational
  Mathematics and Mathematical Physics 7~(1) (1967) 79--120.

\bibitem{selvadurai1979elastic}
A.~P.~S. Selvadurai, Elastic Analysis of Soil-Foundation Interaction, Elsevier
  Science, 1979.

\bibitem{SamarskiiNikolaev1978}
A.~A. Samarskii, E.~S. Nikolaev, Numerical methods for grid equations. Vol. I,
  II, Birkhauser Verlag, Basel, 1989.

\end{thebibliography}
\end{document}